\documentclass[11pt]{amsart}

\usepackage{array,float}
\usepackage{amssymb, amsmath, amsthm, amsbsy, amscd}

\theoremstyle{plain}

\newtheorem{theorem}{Theorem}
\newtheorem*{corollary}{Corollary}

\theoremstyle{remark}
\newcommand{\N}{\ensuremath{\mathbb{N}}}
\newcommand{\Zp}{\ensuremath{\mathbb{Z}_p}}
\newcommand{\Qp}{\ensuremath{\mathbb{Q}_p}}

\newcommand{\Z}{\ensuremath{\mathbb{Z}}}

\newcommand{\mfG}{\ensuremath{\mathfrak{G}}}

\newcommand{\mcN}{{\mathcal{N}}}

\DeclareMathOperator{\Stab}{Stab}

\DeclareMathOperator{\Mat}{Mat}

\DeclareMathOperator{\Aut}{Aut}

\DeclareMathOperator{\G}{\Gamma}
\DeclareMathOperator{\GL}{GL}

\renewcommand{\epsilon}{\varepsilon}
\renewcommand{\phi}{\varphi}

\begin{document}

\title{Enumerating finite class-$2$-nilpotent groups on $2$ generators}

\author{Christopher Voll}

\address{Christopher Voll, School of Mathematics, University of
Southampton, University Road, Southampton SO17 1BJ, United Kingdom}
\email{C.Voll.98@cantab.net}



\begin{abstract}
We compute the numbers $g(n,2,2)$ of nilpotent groups of order
$n$, of class at most $2$ generated by at most $2$ generators, by
giving an explicit formula for the Dirichlet generating function
$\sum_{n=1}^\infty g(n,2,2)n^{-s}$. 



\end{abstract}



\maketitle

\section{Introduction and statement of results}

In \cite{duS/02}, du Sautoy shows how G. Higman's PORC-conjecture
on the numbers $f(n,p)$ of isomorphism types of $p$-groups of
order $p^n$ can be studied using various Dirichlet generating
functions (or zeta functions) associated with groups. Higman
conjectured that the numbers $f(n,p)$ should be `polynomial on
residue classes', i.e.\ that for all $n$ there should exist an
integer $N=N(n)$ and polynomials
$f_1(X),\dots,f_N(X)\in\mathbb{Z}[X]$ such that $f(n,p)=f_i(p)$ if
$p\equiv i$ modulo~$N$. For positive integers $c$ and $d$, we
define $\mcN(c,d)$ to be the set of finite nilpotent groups (up to
isomorphism) of class at most $c$ generated by at most $d$
generators, and put $$g(n,c,d):=\#\{G\in \mcN(c,d)|\;|G|=n\}.$$ We
define the Dirichlet
  generating function
$$\zeta_{c,d}(s):=\sum_{n=1}^\infty g(n,c,d) n^{-s},$$ where $s$
is a
  complex variable. In~\cite[Theorem 1.5]{duS/02} du Sautoy shows that
\begin{equation}\label{euler}
\zeta_{c,d}(s)=\prod_{p \textrm{ prime}}\zeta_{c,d,p}(s),
\end{equation}
where, for a prime $p$, $$\zeta_{c,d,p}(s):=\sum_{i=0}^\infty
g(p^i,c,d) p^{-is}.$$ This `Euler product' reflects the fact that
finite nilpotent groups are the direct products of their Sylow
$p$-subgroups. Du Sautoy goes on to prove that, for all primes~$p$,
the function $\zeta_{c,d,p}(s)$ is rational in~$p^{-s}$ (\cite[Theorem
1.6]{duS/02}). In the remainder of \cite[Part I]{duS/02} he shows that
these local zeta functions are amenable to methods from model theory
and algebraic geometry, and explains how in this setup Higman's
conjecture translates into a question about the reduction modulo $p$
of various ($\Z$-defined) algebraic varieties.

Almost none of the functions $\zeta_{c,d}(s)$ have been explicitly
calculated so far. For $c=1$, it is a trivial and well-known
consequence of the structure theorem for finitely generated
abelian groups that, for all $d\in\N$,
$$\zeta_{1,d}(s)=\prod_{i=1}^d\zeta(is),$$ where
$\zeta(s)=\sum_{n=1}^\infty n^{-s}$ is the Riemann zeta function
(\cite[p. 65]{duS/02}). The purpose of the current note is to give
a formula for the generating function $\zeta_{2,2}(s)$. This seems
to be the only one among the functions $\zeta_{c,d}(s)$, $c>1$,
for which explicit computations exist (see also~\cite[Problem
4]{duS/02}).

\begin{theorem}\label{main theorem} For $(c,d)=(2,2)$ we have $$\zeta_{2,2}(s)=\zeta(s)\zeta(2s)\zeta(3s)^2\zeta(4s).$$
\end{theorem}

\begin{corollary}
For all primes $p$, we have
\begin{multline*}
\zeta_{2,2,p}(s)=1+t+2t^2+4t^3+6t^4+8t^5+13t^6+17t^7+23t^8+31t^9+40t^{10}+\\50t^{11}+65t^{12}+79t^{13}+97t^{14}+119t^{15}+143t^{16}+169t^{17}+203t^{18}+237t^{19}+ O(t^{20})
\end{multline*}
(where $t=p^{-s}$). The abscissa of convergence of $\zeta_{2,2}(s)$ is
$\alpha=1$, and $\zeta_{2,2}(s)$ has an analytic continuation to the
whole complex plane. This continued function has a simple pole at
$s=1$, and thus
$$\sum_{m=1}^ng(m,2,2)\sim\frac{\pi^6}{540}\zeta(3)^2n=2.5725\dots
n.$$ (Here $f(n)\sim g(n)$ means that
$\lim_{n\rightarrow\infty}f(n)/g(n)=1$.)
\end{corollary}

We note that -- as in the case $c=1$ -- the local functions
$\zeta_{2,2,p}(s)$ are rational in $p^{-s}$ with \emph{constant}
coefficients (as the prime $p$ varies). We expect this `strong
uniformity' in the prime $p$ to be the exception rather than the
rule. Indeed, computer calculations of $g(p^i,2,3)$ for
$p\in\{2,3,5\}$ and small values of $i$ show that these numbers do
depend on the prime.\footnote{We thank Eamonn O'Brien for pointing
this out to us.} A `restricted' analogue of Higman's PORC-conjecture
would ask whether, for each pair $(c,d)$, the coefficients of the
functions $\zeta_{c,d,p}(s)$ as rational functions in $p^{-s}$ are
polynomial in $p$ on residue classes of $p$ modulo $N=N(c,d)$. We note
that it is well-known (\cite[Theorem~2]{GSS/88}) that, for all $d\geq
2$ and all primes $p$, the local zeta functions
$\zeta^\triangleleft_{F_{2,d},p}(s)$ enumerating normal subgroups of
finite $p$-power index in the free nilpotent groups $F_{2,d}$ are
rational functions in $p$ and $p^{-s}$ (see also~\cite{Voll/05a}). We
refer to~\cite{duS/02} for an explanation of a link between these
Dirichlet generating functions and the functions $\zeta_{2,d}(s)$
defined above.

In our proof of Theorem~\ref{main theorem} we adopt the strategy
outlined in~\cite[Section 2]{duS/02}. It proceeds by finding
`normal form'-representatives of certain double cosets of integral
matrices, and avoids any algebro-geometric or model-theoretic
considerations. We do, however, use methods pioneered
in~\cite{GSS/88}.

\section{Proof of Theorem~\ref{main theorem}}
The `Euler product decomposition' \eqref{euler} above reduces the
problem to the enumeration of (isomorphism classes of) $p$-groups
of class at most $2$ generated by at most $2$ generators, for a
fixed prime~$p$. Each of these groups occurs as the quotient of
the group $\widehat{F}_p$, the pro-$p$-completion of the free
nilpotent group $F:=F_{2,2}$ of class~$2$ on $2$ generators, by a
normal subgroup $N$ of finite index. The automorphism group
$\mathfrak{G}_p$ of $\widehat{F}_p$ acts on the lattice of these
normal subgroups, and it is known that two subgroups give rise to
the same (isomorphism type of) quotient if and only if they are
equivalent to each other under this action (\cite[Proposition
2.5]{duS/02}). To summarize, we have that, for all primes $p$,
$$\zeta_{2,2,p}(s)=\sum_{N\triangleleft
  \widehat{F}_p}|\widehat{F}_p:N|^{-s}|\mfG_p:\Stab_{\mfG_p}(N)|^{-1}$$
  (cf.~\cite[Theorem 1.13]{duS/02}). To turn this into an explicit
  formula for $\zeta_{2,2,p}(s)$, we firstly linearize the problem of
  counting normal subgroups in the group $\widehat{F}_p$ up to the
  action of~$\mfG_p$, in the following way. Consider the `Heisenberg
  Lie ring' $L:=L_{2,2}:=F/Z(F)\oplus Z(F)$ associated with $F$, with
  Lie bracket induced from taking commutators in $F$. 
It is well-known (\cite[Remark on p. 206]{GSS/88}) that normal
subgroups of index $p^n$ in $\widehat{F}_p$ correspond to ideals of
index $p^n$ in the $\Zp$-Lie algebra $L_p:=L\otimes\Zp$, and it is
easily verified that orbits under $\mfG_p$ in the lattice of normal
subgroups of $\widehat{F}_p$ correspond to orbits under
$\Aut(L_p)$. Having chosen a basis for $L$, e.g.\ $(x,y,z)$ with
$[x,y]=z$, full sublattices $\Lambda$ in $(L_p,+)$ may be identified
with cosets $\G M$, $\G:=\GL_3(\Zp)$, $M\in\Mat(3,\Zp)\cap\GL_3(\Qp)$,
by encoding in the rows of $M$ the coordinates (with respect to the
chosen basis, viewed as a basis for the $\Zp$-Lie algebra $L_p$) of
generators for $\Lambda$.

Given the basis $(x,y,z)$ as above, every coset $\G M$
corresponding to a lattice of finite index in $L_p$ contains a
unique matrix of the form
\begin{equation}\label{reduced form}
M=\left( \begin{array}{ccc} p^{n_1}&a_{12}&a_{13} \\
&p^{n_2}&a_{23}
                      \\ & &p^{n_3}
            \end{array} \right),
\end{equation}
with $n_1,n_2,n_3\in\N_0:=\{0,1,2,\dots\}$, $0\leq a_{12}<
p^{n_2}$ and~$0 \leq a_{13}, a_{23}< p^{n_3}$. It is easy to see
and well-known that $\G M$ corresponds to an {\sl ideal} in $L_p$
(of index $p^{n_1+n_2+n_3}$) if and only if
\begin{equation}\label{ideal condition}
n_3 \leq {n_2},{n_1}, v_p(a_{12})
\end{equation}
(where $v_p$ denotes the $p$-adic valuation).

The choice of basis allows us to identify $\Aut(L_p)$ with the
group of matrices of the form
$$\left\{\left(\begin{matrix}\alpha&*\\&\det(\alpha)\end{matrix}\right)|\;\alpha\in\GL_2(\Zp)\right\}\subseteq\G.$$
The action of $\Aut(L_p)$ on the lattice of ideals of finite index in
$L_p$ is then simply given by the natural (right-)action on the set of
cosets $\G M$. We claim that the matrices of the form \eqref{reduced
form} with
\begin{equation}\label{normal form}
n_1=e_1+e_2+e_3,\quad n_2=e_2+e_3, n_3=e_3, \quad a_{12}=0, \quad a_{13}=p^{e_4}, \quad a_{23}=p^{e_5},
\end{equation}
 where
\begin{equation}
e_1,\dots,e_5 \in\N_0, \quad {e_4},{e_5} \leq {e_3},\quad e_5\leq e_4
\leq e_5+e_1\label{first condition}
\end{equation}
form a complete set of representatives of the double cosets
$$\G\backslash (\Mat(3,\Zp)\cap\GL_3(\Qp)) /\Aut(L_p).$$ This
suffices as then
\begin{align*}
\zeta_{2,2,p}(s)&=\sum_{\G\backslash M
/\Aut(L_p)}|\det(M)|^{-s}\\&=\sum_{(e_1,\dots,e_5)\in \N_0^5\textrm{
satisfying }\eqref{first
condition}}p^{-3{e_3}s-2{e_2}s-e_1s}\\&=\frac{1}{(1-p^{-s})(1-p^{-2s})(1-p^{-3s})^2(1-p^{-4s})}\\&=\zeta_p{(s)}\zeta_p{(2s)}\zeta_p{(3s)}^2\zeta_p{(4s)},
\end{align*}
as an easy calculation yields.

To prove our claim, it suffices to show that every double coset
contains a matrix of the form \eqref{reduced form},
satisfying~\eqref{normal form} and \eqref{first condition}, and that,
if $N$ and $N'$ are in this normal form, the double cosets defined by
them coincide only if $N=N'$.

We start by observing that every double coset of matrices of the
form~\eqref{reduced form} satisfying the `ideal
condition'~\eqref{ideal condition} contains a matrix of the form~\eqref{reduced form} with
$$n_1 = e_1 + e_2 + e_3, \quad n_2 = e_2 + e_3, \quad n_3=e_3, \quad
a_{12}=0,$$ where $$e_1,e_2,e_3\in\N_0,\quad e_4:=v_p(a_{13})\leq
e_3,\quad e_5:=v_p(a_{23})\leq e_3.$$ This is because $\Aut(L_p)$
contains a copy of $\GL_2(\Zp)$, allowing us to bring the top-left
$2\times2$-block of $N$ into `Smith normal form'.
By a suitable base change we may also arrange for $a_{13}=p^{e_4},
a_{23}=p^{e_5}$. Furthermore we can achieve that $e_5\leq e_4 \leq
e_5+e_1$ as we can always add multiples of the first row to the second
row and multiples of $p^{e_1}$ times the first row to the second row,
each time `clearing' the $a_{21}$- and $a_{12}$-entry respectively by
right-multiplication by suitable elements in $\Aut(L_p)$, namely
elementary column operations involving only the first two columns,
leaving the third column stable.

Now assume that we are given matrices $N$ and $N'$ in normal
form~\eqref{reduced form}, satisfying \eqref{normal form} and
\eqref{first condition}, with associated invariants $(e_1,\dots,e_5)$
and $(e_1',\dots,e_5')$, respectively. For $N$ and $N'$ to define the
same double cosets it is clearly necessary that
$(e_1,e_2,e_3)=(e_1',e_2',e_3')$. Also, necessarily $e_5=e_5'$, as
$e_5$ determines the $p$-adic norm of the last column of $N$, which is
invariant both by left-multiplication by elements in $\G$ and by
right-multiplication by $\Aut(L_p)$. We thus have, without loss of
generality, 
\begin{gather*}N=\left(
\begin{array}{ccc} p^{{e_1}+{e_2}+{e_3}}& &p^{e_5+m+n} \\
                      &p^{{e_2}+{e_3}}&p^{e_5} \\ & &p^{e_3}
            \end{array} \right), \\
 N'=\left( \begin{array}{ccc}
                      p^{{e_1}+{e_2}+{e_3}}&               &p^{e_5+m}   \\
                                           &p^{{e_2}+{e_3}}&p^{e_5}   \\
                                           &               &p^{e_3}
            \end{array} \right),
\end{gather*}
where $m,n\geq 0,\; m+n \leq e_3$ and therefore in particular
$m<e_3$. If $N$ and $N'$ were to define the same double coset, we must
be able to `adjust' the $a_{13}$-entry of $N$ to have valuation
$p^{e_5+m}$ by a suitable row operation. The only way to achieve this
is to add $up^m$ times the second row to the first, where $u$ is a
$p$-adic unit.
If $m>0$, there is no other way to achieve Smith normal form on the
top-left $2\times2$-block than to reverse this row operation. If $m=0$
the only alternative is to choose the $a_{12}$-entry as a pivot to
obtain this, 
which again has the invariants of $N$, not of $N'$ if $n\neq 0$.
Thus $n=0$, completing the proof of Theorem~\ref{main theorem}.


\section*{Acknowledgements}
The result presented in this note forms part of my Cambridge Ph.D.
thesis~\cite{Voll/02}, which was supervised by Marcus du Sautoy
and partly funded by the Studienstiftung des deutschen Volkes and
the Cambridge European Trust. I would like to thank Eamonn O'Brien
for providing me with various computer-generated numerical data,
which proved immensely helpful. These computations are based
on~\cite{OBrien/90}.

\bibliographystyle{amsplain}
\def\cprime{$'$}
\providecommand{\bysame}{\leavevmode\hbox to3em{\hrulefill}\thinspace}
\providecommand{\MR}{\relax\ifhmode\unskip\space\fi MR }
\providecommand{\MRhref}[2]{%
  \href{http://www.ams.org/mathscinet-getitem?mr=#1}{#2}
}
\providecommand{\href}[2]{#2}

\end{document}